\documentclass[preprint,1p,times]{elsarticle}
\usepackage{amsmath, amssymb, amsthm}
\usepackage{graphicx}
\usepackage{color,soul}
\usepackage{float}
\usepackage{lineno}
\usepackage{pgf}
\usepackage{epstopdf}

\journal{ }
\begin{document}

\begin{frontmatter}
	
\title{A fractional calculus approach to Rosenzweig-MacArthur predator-prey model and its solution}
\author[cmbe]{Shuvojit Mondal\fnref{fn1}}
\ead{shuvojitmondal91@gmail.com}
\author[cmbe]{Nandadulal Bairagi\fnref{fn2}}
\ead{nbairagi.math@jadavpuruniversity.in}
\author[ju]{Abhijit Lahiri\corref{cor1}}
\ead{lahiriabhijit2000@yahoo.com}
\cortext[cor1]{Corresponding author}
%\fnteyt[fn2]{Research is supported by PURSE DST Phase II.}
\address[cmbe]{Centre for Mathematical Biology and Ecology\\ Department of Mathematics, Jadavpur University\\ Kolkata-700032, India.}
\address[ju]{Department of Mathematics, Jadavpur University \\ Kolkata-700032, India.}

\begin{abstract}
In this paper we present analytical solution of a fractional order predator-prey model, where prey grows logistically and predation occurs following type II response function, by homotopy perturbation method. Numerical solutions are presented to illustrate different particular cases. Our computational results show that accurate solution may be obtained with few iterations.
\end{abstract}

\begin{keyword}
Predator-prey model, Caputo fractional derivative, Homotopy perturbation method, Initial value problem.
\end{keyword}

\end{frontmatter}

\section{Introduction}
%There are several kind of approaches of fractional calculus (differentiation and integration) to the various ecological models which are done by many researchers in recent years.
Fractional calculus, which is generalization of integer order differentiation and $n-$fold integration, has been successfully applied in different branches of science and engineering. In recent past, it has also been applied in several ecological models \cite{CuiYang14, DasGupta09, DasGupta11}. Differential equations with fractional-order derivatives (or integrals) are generally called fractional differential (or integral) equations. \\
If $x(t)$ and $y(t)$ be, respectively, the densities of prey and predator populations at time $t$ then a general model for predator-prey interaction can be written as
\begin{eqnarray}\label{Dynamical_system}\nonumber
\frac{dx}{dt}=xg(x,y), \\
\frac{dy}{dt}=yh(x,y),
\end{eqnarray}
where $g$ and $h$ are the per capita growth rates of prey and predator populations. If we choose $g(x,y) = a-by$ and $h(x,y) = cx-d$, then the above predator-prey system (\ref{Dynamical_system}) becomes the familiar Lotka-Volterra system:
% " The logistic equation" which was presented by Lotka\cite{Lotka25} and Volterra\cite{Volterra26}.
\begin{eqnarray}\label{Lotka-Volterra model}\nonumber
\frac{dx}{dt} & = & ax(t) - bx(t)y(t), \\
\frac{dy}{dt} & = & c x(t)y(t) - dy(t).
\end{eqnarray}
Here prey population grows exponentially with rate constant $a$, predator consumes prey at a rate $b$. The parameters $c$ and $d$ represent, respectively, the growth rate and death rate of predator. May \cite{May72} studied the system (\ref{Lotka-Volterra model}) and showed that the interior equilibrium has neutral stability.
Fractional differential equations have the ability of providing an exact or an approximate solutions of a nonlinear system. The main advantage of fractional order system is that they allow greater degrees of freedom than an integer order system \cite{DasGupta11}. Any first order system of differential equations can be transformed into a fractional order system by replacing just the ordinary derivative by fractional order derivative. Das et al. \cite{DasGupta11} have recently solved the following fractional order Lotka-Volterra predator-prey model by Homotopy perturbation method
\begin{eqnarray}\label{Lotka-Volterra predator-prey model}
D^{\alpha}_{t}x & = & a(t)x(t) - b(t)x(t)y(t), \\
D^{\beta}_{t}y & = & c(t) x(t)y(t) - d(t)y(t),\nonumber
\end{eqnarray}
with fractional orders $\alpha$ and $\beta$, where $0<\alpha\leq1$, $0<\beta\leq1$ and $D^{\alpha}_{t}\equiv \frac{d^\alpha}{dt^\alpha}$ is the fractional derivative of order $\alpha$ in Caputo sense. Here $a, b, c$ and $d$ have been considered as function of time $t$. In an earlier study \cite{DasGupta09}, they studied the same model by Homotopy perturbation method assuming that $a=b=c=d=1$. Cui and Yang \cite{CuiYang14} modified the system (\ref{Lotka-Volterra predator-prey model}) by considering predator's density dependent death.
All these models assume that the predator's functional response, defined as the number of prey catch per predator per unit of time (here it is $bx$), is unsaturated. It implies that per capita predation increases with the prey density and there exists no upper limit to the prey consumption rate of predator. In natural system, however, per capita prey consumption should satiate as the prey density increases. Holling type II response function of predator represented by $\frac{\alpha x}{a+x}$, where $\alpha$ is the saturation vale of the functional response and $a$ is the half-saturation constant, is assumed to be a more appropriate function to express predator's prey capture rate \cite{BET07}. To the best of our knowledge, nobody has considered this type II function to find the analytical solution of a fractional order predator-prey model as the analytical complexity multiplies in this case. Here we find the analytical solution of a predator-prey model when predator's response function is type II. We consider the following benchmark integer order predator-prey model, popularly known as Rosenzweig-MacArthur model \cite{SS12}:
\begin{eqnarray}\label{Rosenzweig-MacArthur model}
\frac{dx}{dt} & = & rx\bigg(1-\frac{x}{K}\bigg) - \frac{\alpha x y}{a+x}, \\
\frac{dy}{dt}& = & \frac{\beta x y}{a+x} - dy. \nonumber
\end{eqnarray}
This model says that the prey population $x(t)$ grows in logistic fashion with intrinsic growth rate $r$ and carrying capacity $K$. Predator $y(t)$ captures the prey at a maximum rate $\alpha$ with half saturation constant $a$. The parameters $\beta$ and $d$ represent, respectively, the growth and death rates of predator. \\
Considering the fractional derivatives in the sense of Caputo derivative and assuming $0<m\leq1, 0<n\leq1$, we have the following fractional order predator-prey model:
\begin{eqnarray}\label{Rosenzweig-MacArthur predator-prey model}
D^{m}_{t}x & = & rx\bigg(1-\frac{x}{K}\bigg) - \frac{\alpha x y}{a+x}, \\
D^{n}_{t}y & = & \frac{\beta x y}{a+x} - dy. \nonumber
\end{eqnarray}
The initial conditions are considered as $x(0) = \delta$ and $y(o) = \gamma$. From biological point of view, $\delta$ and $\gamma$ are assumed to be positive and all other parameters are also considered to be positive. The main advantage of Caputo's approach is that the initial conditions for the fractional differential equations with Caputo derivatives takes the similar form as for integer-order differential equations \cite{CuiYang14, Podlubny99}, i.e., it has advantage of defining integer order initial conditions for fractional order differential equations, but unfortunately Riemann-Liouville fractional derivative approach is unable to define such thing. In this study, we analyze the system (\ref{Rosenzweig-MacArthur predator-prey model}) and find its approximate analytical solution with the help of Homotopy perturbation method.\\

The paper is organized as follows: In Section 2, we give some definitions and basic conception of Homotopy perturbation method. Analytical solution of the model is presented in Section 3. Extensive numerical computation are presented in Section 4 and the paper ends with a summary in Section 5.\\
\section{Preliminaries}
\subsection{Fractional calculus}
Fractional calculus is a generalization of ordinary differentiation and integration to an arbitrary (non-integer) order and therefore justifies its alternative name as calculus of arbitrary order. During last two decades there are several kind of researchers gave the ideas to define fractional calculus, for example  Mathematician Liouville, Riemannn, Grunwald-Letnikov, Caputo have done major work on fractional calculus. we are going to talk about the definitions of fractional integral$/$derivative and some preliminary results in this section. \\

\noindent\textbf{Definition 2.1.} \cite{Podlubny99, MillerRoss93} The fractional integral of order $\nu\geq0$ of a function $f(t):R^+ \longrightarrow R$ is given by \\
$$J_t^\nu f(t) = \frac{1}{\Gamma{(\nu)}} \int_0^t(t-s)^{\nu - 1}f(s)ds,$$ \\
provided the right hand side integral is point-wise continuous on $R^+$. Here $R$ is the set of real numbers, $R^+$ is the set of positive real numbers and $\Gamma{(\cdot)}$ is the well-known gamma function. \\

\noindent\textbf{Definition 2.2.} \cite{Podlubny99, MillerRoss93} The caputo derivative of order $\mu>0$ of a continuous  function $\phi(t) :R^+ \longrightarrow R$ is defined by \\
$$D_t^\mu \phi(t) = \frac{1}{\Gamma{(n-\mu)}} \int_0^t\frac{\phi^{(n)}(\tau)}{(t-\tau)^{\mu-n+1}}d\tau,$$ \\
provided the right hand side integral is point-wise continuous on $R^+$. Here $n$ is an integer, $\mu$ is a real number and $n=[\mu]+1$, i.e., $(n-1)\leq\mu<n$.\\

Note that $J_t^\nu$ is the integral operator of order $\nu>0$, whereas $D_t^\mu$ is the differential operator, specifically the inverse of integral operator. Various kind of properties of these integral and differential operators can be found in \cite{MillerRoss93, Podlubny99} and we reproduce here some of them.\\
For $f(t)\in C_\theta, \theta\geq-1, \rho, \sigma\geq0$ and $\lambda>-1$, we have the following properties:
\begin{itemize}
	\item[(1)]$ J^\rho J^\sigma f(t) = J^{\rho+\sigma} f(t)$,
	\item[(2)] $J^\rho J^\sigma f(t) = J^\sigma J^\rho f(t)$,
	\item[(3)] $J^\rho t^\lambda  = \frac{\Gamma{(\lambda+1)}}{\Gamma{(\rho+\lambda+1)}} t^{\rho+\lambda}$,
	\item[(4)] $D^\sigma J^\rho f(t)  =  J^{\rho-\sigma} f(t)$. \end{itemize}
\subsection{Homotopy Perturbation Method}
Homotopy perturbation method (HPM), first proposed by He \cite{HE99,HE00}, has been employed to solve a large variety of linear and non linear problems \cite{HE00,Taghipour11,HE04,HE05,HE98,Ganji06}. Generally, this method is a coupling of normal perturbation method and homotopy in topology. Contrary to other method, like Adomian decomposition method (ADM), this method does not require a small parameter in the equation and considered as the main advantage to get the analytical approximate solution easily and elegantly without transforming the equation or linearizing the problem. HPM gives a very rapid convergence of the series solution, generally after a few number of iterations, and leads to approximate solutions similar to the accurate solutions of nonlinear problems \cite{MomaniOdibat07}. We, therefore, apply HPM to find the approximate analytical solution of our system (\ref{Rosenzweig-MacArthur predator-prey model}).\\
To discus this method, we consider the following nonlinear differential equation
\begin{equation}\label{Homotopy perturbation method_1}
L(u)+N(u) = f(r), r\in\Omega
\end{equation}
with the boundary conditions
\begin{equation}\label{Homotopy perturbation method_2}
B(u,\frac{\partial u}{\partial n}) = 0, r\in\Gamma,
\end{equation}
where $L$ is the linear operator, $N$ is the nonlinear operator, $B$ is the boundary operator, $\Gamma$ is the boundary of the domain $\Omega$ and $f(r)$ is the known function.\\
Following He's homotopy perturbation technique \cite{HE99,HE00}, we construct a homotopy as $v(r,p):\Omega\times[0,1]\rightarrow R$ which satisfies\\
\begin{equation}\label{Homotopy perturbation method_3}
H(v,p) = (1-p)[L(v)-L(u_0)] + p[L(v)+N(v)-f(r)]\\
= 0, p\in[0,1], r\in\Omega
\end{equation}
or\\
\begin{equation}\label{Homotopy perturbation method_4}
H(v,p) = [L(v)-L(u_0)] + pL(u_0)+p[N(v)-f(r)]\\
= 0, p\in[0,1], r\in\Omega.
\end{equation}
Here, $p\in[0,1]$ is the embedding parameter and $u_0$ is the initial approximation which satisfies the boundary conditions. From (\ref{Homotopy perturbation method_3}) and (\ref{Homotopy perturbation method_4}), one can have
\begin{equation}\label{Homotopy perturbation method_5}
H(v,0) = L(v)-L(u_0) = 0,
\end{equation}
\begin{equation}\label{Homotopy perturbation method_6}
H(v,1) = L(v)+N(v)-f(r) = 0.
\end{equation}
It means that the changing process of the embedding parameter $p$ from zero to unity is just that of $v(r,p)$ from $u_0$ to $u(r)$. From the topological point of view this is called deformation and $L(v)-L(u_0)$, $L(v)+N(v)-f(r)$ are called homotopic.\\
According to HPM \cite{MomaniOdibat07, Taghipour11}, we assume that the embedding parameter $p$ $(0\leq p\leq 1)$  as a "small parameter" and also assume that the solutions of the Eqs. (\ref{Homotopy perturbation method_3}) and (\ref{Homotopy perturbation method_4}) can be expressed as a power series in $p$:
\begin{equation}\label{Homotopy perturbation method_7}
v = \sum_{n = 0}^\infty p^n v_n(t) = v_0(t) + pv_1(t) + p^2v_2(t) + p^3v_3(t) +.......\\
\end{equation}
Now the approximate solution of Eqs. (\ref{Homotopy perturbation method_1}) can be easily obtained by setting $p\rightarrow 1$:
\begin{equation}\label{Homotopy perturbation method_8}
u = \lim_{p\rightarrow1}v = \lim_{p\rightarrow1}\sum_{n = 0}^\infty p^n v_n(t) = \sum_{n = 0}^\infty v_n(t) = v_0(t) + v_1(t) + v_2(t) + v_3(t) +.......\\
\end{equation}
The convergence of the series in (\ref{Homotopy perturbation method_8}) has been proved in He's paper \cite{HE99,HE00}. The combination of the perturbation method and the homotpoy method is called the homotpoy perturbation method. Though the convergence depends on the nonlinear operator $N$ in (\ref{Homotopy perturbation method_1}), the series in (\ref{Homotopy perturbation method_8}) is convergent in most cases.
\section{Analytical solution of the problem}
In this section, we apply the Homotopoy perturbation method to solve the Rosenzweig and MacArthur predator-prey model (\ref{Rosenzweig-MacArthur predator-prey model}) with initial conditions $x(0) = \delta$ and $y(o) = \gamma$.\\
Following HPM,  we construct the homotopy structure  of (\ref{Rosenzweig-MacArthur predator-prey model}) as follows:
\begin{eqnarray}\label{Rosenzweig-MacArthur predator-prey model_1}\nonumber
D^{m}_{t}x & = & p\bigg[rx\bigg(1-\frac{x}{K}\bigg) - \frac{\alpha x y}{a+x}\bigg], \\
D^{n}_{t}y & = & p\bigg[\frac{\beta x y}{a+x} - dy\bigg], \\ \nonumber
\end{eqnarray}
where $0<m,n\leq1$ and $p\in[0,1]$ is the homotopy parameter. If we consider $p = 0$, then (\ref{Rosenzweig-MacArthur predator-prey model_1}) becomes a homogeneous fractional differential equation, which can be easily solved by fractional approach \cite{Podlubny99, MillerRoss93}.  The system (\ref{Rosenzweig-MacArthur predator-prey model_1}) returns to the original Eqn. (\ref{Rosenzweig-MacArthur predator-prey model}) for $p = 1$. \\

Following HPM, we assume that solutions of (\ref{Rosenzweig-MacArthur predator-prey model_1}) can be written as a power series in $p$:
\begin{equation}\label{Rosenzweig-MacArthur predator-prey model_2}
x(t) = \sum_{n = 0}^\infty p^n x_n(t) =   x_0(t) + px_1(t) + p^2x_2(t) + p^3x_3(t) +......., \\
\end{equation}
\begin{equation}\label{Rosenzweig-MacArthur predator-prey model_3}
y(t) = \sum_{n = 0}^\infty p^n y_n(t) =  y_0(t) + py_1(t) + p^2y_2(t) + p^3y_3(t) +......., \\
\end{equation}
Taking $p \rightarrow 1$, we obtain the approximate solutions of the original Eqn.(\ref{Rosenzweig-MacArthur predator-prey model}) as
\begin{equation}\label{Rosenzweig-MacArthur predator-prey model_5}
x(t) = \lim_{p\rightarrow1}\sum_{n = 0}^\infty p^n x_n(t) =   x_0(t) + x_1(t) + x_2(t) + x_3(t) +......., \\
\end{equation}
\begin{equation}\label{Rosenzweig-MacArthur predator-prey model_6}
y(t) = \lim_{p\rightarrow1} \sum_{n = 0}^\infty p^n y_n(t) =  y_0(t) + y_1(t) + y_2(t) + y_3(t) +......., \\
\end{equation}

Substituting Eqs. (\ref{Rosenzweig-MacArthur predator-prey model_2}), (\ref{Rosenzweig-MacArthur predator-prey model_3}) in (\ref{Rosenzweig-MacArthur predator-prey model_1}) and equating the powers of $p$ from both sides, we obtain the following set of linear fractional order differential equations:
\begin{eqnarray}\label{Rosenzweig-MacArthur predator-prey model_4}
\begin{split}
p^0 : D^m x_0(t) =&0, \\
D^n y_0(t) =&0, \\
p^1 : D^m x_1(t) =&\bigg(rx_0-\frac{rx_0^2}{K}\bigg)-\alpha x_0 y_0\bigg(\frac{1}{a}-\frac{x_0}{a^2}+\frac{x_0^2}{a^3}-\frac{x_0^3}{a^4}\bigg), \\
D^n y_1(t) =&\beta x_0 y_0 \bigg(\frac{1}{a}-\frac{x_0}{a^2}+\frac{x_0^2}{a^3}-\frac{x_0^3}{a^4}\bigg)-d y_0, \\
p^2 : D^m x_2(t) =&\bigg(rx_1-\frac{2rx_0x_1}{K}\bigg)-\alpha x_0y_0\bigg(\frac{2x_0x_1}{a^3}-\frac{x_1}{a^2}-\frac{3x_0^2x_1}{a^4}\bigg)\\
&-\alpha(x_0y_1+x_1y_0)\bigg(\frac{1}{a}-
\frac{x_0}{a^2}+\frac{x_0^2}{a^3}-\frac{x_0^3}{a^4}\bigg), \\
D^n y_2(t) =&\beta x_0 y_0\bigg(\frac{2x_0x_1}{a^3}-\frac{x_1}{a^2}-\frac{3x_0^2x_1}{a^4}\bigg)\\
&+\beta (x_0y_1+x_1y_0)\bigg(\frac{1}{a}-
\frac{x_0}{a^2}+\frac{x_0^2}{a^3}-\frac{x_0^3}{a^4}\bigg)-dy_1, \\
p^3 : D^m x_3(t) =&\bigg(rx_2-\frac{r}{K}(x_1^2+2x_0x_2)\bigg)-\alpha x_0y_0\bigg(\frac{x_1^2+2x_0x_2}{a^3}-\frac{x_2}{a^2}-\frac{3(x_0^2x_2+x_0x_1^2)}{a^4}\bigg) \\
&-\alpha(x_0y_1+x_1y_0)\bigg(\frac{1}{a}-
\frac{x_0}{a^2}+\frac{x_0^2}{a^3}-\frac{x_0^3}{a^4}\bigg)\\
&-\alpha(x_0y_2+x_2y_0+x_1y_1)\bigg(\frac{1}{a}-
\frac{x_0}{a^2}+\frac{x_0^2}{a^3}-\frac{x_0^3}{a^4}\bigg), \\
D^n y_3(t) =&\beta x_0y_0\bigg(\frac{x_1^2+2x_0x_2}{a^3}-\frac{x_2}{a^2}-\frac{3(x_0^2x_2+x_0x_1^2)}{a^4}\bigg)\\
&+\beta(x_0y_1+x_1y_0)\bigg(\frac{1}{a}-
\frac{x_0}{a^2}+\frac{x_0^2}{a^3}-\frac{x_0^3}{a^4}\bigg)\\
&+\beta(x_0y_2+x_2y_0+x_1y_1)\bigg(\frac{1}{a}-
\frac{x_0}{a^2}+\frac{x_0^2}{a^3}-\frac{x_0^3}{a^4}\bigg)-dy_2, \\
\end{split}
\end{eqnarray}
and so on.\\
We now apply the operators $J^m_t$ and $J^n_t$, which are basically represented as the inverse operators of the Caputo derivative $D^m_t$ and $D^n_t$ respectively, on both sides of each fractional differential equations (\ref{Rosenzweig-MacArthur predator-prey model_4}).
%From the theory of homotopy perturbation method, the approximate solutions for $x(t)$ and $y(t)$ by setting $p \rightarrow 1$ as
%\begin{equation}\label{Rosenzweig-MacArthur predator-prey model_7}
%x(t) = \lim_{N_0\rightarrow\infty} \sum_{n=0}^{N_0-1} x_n(t), \\
%\end{equation}
%\begin{equation}\label{Rosenzweig-MacArthur predator-prey model_8}
%y(t) = \lim_{N_0\rightarrow\infty} \sum_{n=0}^{N_0-1} y_n(t), \\
%\end{equation}
%where $N_0\geq1$.
%%the main thing is that the above two series converge very rapidly in real physical problems. The rapid convergence means that only few terms are required to get approximate solution, so therefore there requires no infinite series.\\
Solving each equation of (\ref{Rosenzweig-MacArthur predator-prey model_4}), we have
\begin{equation}\nonumber
\begin{split}
x_0(t) =& \delta, \\
y_0(t) =& \gamma, \\
x_1(t) =&
\bigg(r\delta-\frac{r\delta^2}{K}-\alpha\delta\gamma\bigg(\frac{1}{a}-\frac{\delta}{a^2}+\frac{\delta^2}{a^3}-\frac{\delta^3}{a^4}\bigg)\bigg)\frac{t^m}{\Gamma(m+1)}, \\
y_1(t) =&
\bigg(\beta\delta\gamma\bigg(\frac{1}{a}-\frac{\delta}{a^2}+\frac{\delta^2}{a^3}-\frac{\delta^3}{a^4}\bigg)-d\gamma\bigg)\frac{t^n}{\Gamma(n+1)}, \\
x_2(t) =&
\bigg(r-\frac{2r\delta}{K}-\alpha\delta\gamma\bigg(\frac{2\delta}{a^3}-\frac{1}{a^2}-\frac{3\delta^2}{a^4}\bigg)-\alpha\gamma\bigg(\frac{1}{a}-\frac{\delta}{a^2}+\frac{\delta^2}{a^3}-\frac{\delta^3}{a^4}\bigg)\bigg)\\
&\bigg(r\delta-\frac{r\delta^2}{K}-\alpha\delta\gamma\bigg(\frac{1}{a}-\frac{\delta}{a^2}+\frac{\delta^2}{a^3}-\frac{\delta^3}{a^4}\bigg)\bigg)\frac{t^{2m}}{\Gamma(2m+1)}\\
&-\alpha\delta\bigg(\frac{1}{a}-\frac{\delta}{a^2}+\frac{\delta^2}{a^3}-\frac{\delta^3}{a^4}\bigg)
\bigg(\beta\delta\gamma\bigg(\frac{1}{a}-\frac{\delta}{a^2}+\frac{\delta^2}{a^3}-\frac{\delta^3}{a^4}\bigg)-d\gamma\bigg)\frac{t^{m+n}}{\Gamma(m+n+1)}, \\
y_2(t) =&
\bigg(\beta\delta\gamma\bigg(\frac{2\delta}{a^3}-\frac{1}{a^2}-\frac{3\delta^2}{a^4}\bigg)+\beta\gamma\bigg(\frac{1}{a}-\frac{\delta}{a^2}+\frac{\delta^2}{a^3}-\frac{\delta^3}{a^4}\bigg)\bigg)\\
&\bigg(r\delta-\frac{r\delta^2}{K}-\alpha\delta\gamma\bigg(\frac{1}{a}-\frac{\delta}{a^2}+\frac{\delta^2}{a^3}-\frac{\delta^3}{a^4}\bigg)\bigg)\frac{t^{m+n}}{\Gamma(m+n+1)}\\
&+\bigg(\beta\delta\bigg(\frac{1}{a}-\frac{\delta}{a^2}+\frac{\delta^2}{a^3}-\frac{\delta^3}{a^4}\bigg)-d\bigg)
\bigg(\beta\delta\gamma\bigg(\frac{1}{a}-\frac{\delta}{a^2}+\frac{\delta^2}{a^3}-\frac{\delta^3}{a^4}\bigg)-d\gamma\bigg)\frac{t^{2n}}{\Gamma(2n+1)}, \\
\end{split}
\end{equation}
\begin{equation}\nonumber
\begin{split}
x_3(t) =&
\bigg(r-\frac{2r\delta}{K}-\alpha\delta\gamma\bigg(\frac{2\delta}{a^3}-\frac{1}{a^2}-\frac{3\delta^2}{a^4}\bigg)-\alpha\gamma\bigg(\frac{1}{a}-\frac{\delta}{a^2}+\frac{\delta^2}{a^3}-\frac{\delta^3}{a^4}\bigg)\bigg)^2\\
&\bigg(r\delta-\frac{r\delta^2}{K}-\alpha\delta\gamma\bigg(\frac{1}{a}-\frac{\delta}{a^2}+\frac{\delta^2}{a^3}-\frac{\delta^3}{a^4}\bigg)\bigg)\frac{t^{3m}}{\Gamma(3m+1)}\\
&-\bigg(\frac{r}{K}+\frac{\alpha\delta\gamma}{a^3}-\frac{3\alpha\delta^2\gamma}{a^4}\bigg)\bigg(r\delta-\frac{r\delta^2}{K}-\alpha\delta\gamma\bigg(\frac{1}{a}-\frac{\delta}{a^2}+\frac{\delta^2}{a^3}-\frac{\delta^3}{a^4}\bigg)\bigg)^2\frac{\Gamma(2m+1)t^{3m}}{\Gamma(m+1)^2\Gamma(3m+1)}\\
&-\alpha\delta\bigg(\frac{1}{a}-\frac{\delta}{a^2}+\frac{\delta^2}{a^3}-\frac{\delta^3}{a^4}\bigg)\bigg(\beta\delta\gamma\bigg(\frac{1}{a}-\frac{\delta}{a^2}+\frac{\delta^2}{a^3}-\frac{\delta^3}{a^4}\bigg)-d\gamma\bigg)\frac{t^{m+n}}{\Gamma(m+n+1)}\\
&-\alpha\gamma\bigg(\frac{1}{a}-\frac{\delta}{a^2}+\frac{\delta^2}{a^3}-\frac{\delta^3}{a^4}\bigg)\bigg(r\delta-\frac{r\delta^2}{K}-\alpha\delta\gamma\bigg(\frac{1}{a}-\frac{\delta}{a^2}+\frac{\delta^2}{a^3}-\frac{\delta^3}{a^4}\bigg)\bigg)\frac{t^{2m}}{\Gamma(2m+1)}\\
&-\alpha\delta\bigg(\frac{1}{a}-\frac{\delta}{a^2}+\frac{\delta^2}{a^3}-\frac{\delta^3}{a^4}\bigg)\bigg(r-\frac{2r\delta}{K}-\alpha\delta\gamma\bigg(\frac{2\delta}{a^3}-\frac{1}{a^2}-\frac{3\delta^2}{a^4}\bigg)-\alpha\gamma\bigg(\frac{1}{a}-\frac{\delta}{a^2}+\frac{\delta^2}{a^3}-\frac{\delta^3}{a^4}\bigg)\bigg)\\
&\bigg(\beta\delta\gamma\bigg(\frac{1}{a}-\frac{\delta}{a^2}+\frac{\delta^2}{a^3}-\frac{\delta^3}{a^4}\bigg)-d\gamma\bigg)\frac{t^{2m+n}}{\Gamma(2m+n+1)}
-\alpha\bigg(\frac{1}{a}-\frac{\delta}{a^2}+\frac{\delta^2}{a^3}-\frac{\delta^3}{a^4}\bigg)\\
&\bigg(r\delta-\frac{r\delta^2}{K}-\alpha\delta\gamma\bigg(\frac{1}{a}-\frac{\delta}{a^2}+\frac{\delta^2}{a^3}-\frac{\delta^3}{a^4}\bigg)\bigg)
\bigg(\beta\delta\gamma\bigg(\frac{1}{a}-\frac{\delta}{a^2}+\frac{\delta^2}{a^3}-\frac{\delta^3}{a^4}\bigg)-d\gamma\bigg)\\
&\frac{\Gamma(m+n+1)t^{2m+n}}{\Gamma(m+1)\Gamma(n+1)\Gamma(2m+n+1)}
\end{split}
\end{equation}
\begin{equation}\nonumber
\begin{split}
&-\alpha\beta\delta\bigg(\frac{1}{a}-\frac{\delta}{a^2}+\frac{\delta^2}{a^3}-\frac{\delta^3}{a^4}\bigg)\bigg(\gamma\bigg(\frac{1}{a}-\frac{\delta}{a^2}+\frac{\delta^2}{a^3}-\frac{\delta^3}{a^4}\bigg)\\
&+\delta\gamma\bigg(\frac{2\delta}{a^3}-\frac{1}{a^2}-\frac{3\delta^2}{a^4}\bigg)\bigg)
\bigg(r\delta-\frac{r\delta^2}{K}-\alpha\delta\gamma\bigg(\frac{1}{a}-\frac{\delta}{a^2}+\frac{\delta^2}{a^3}-\frac{\delta^3}{a^4}\bigg)\bigg)\frac{t^{2m+n}}{\Gamma(2m+n+1)}\\
&-\alpha\delta\bigg(\frac{1}{a}-\frac{\delta}{a^2}+\frac{\delta^2}{a^3}-\frac{\delta^3}{a^4}\bigg)\bigg(\beta\delta\bigg(\frac{1}{a}-\frac{\delta}{a^2}+\frac{\delta^2}{a^3}-\frac{\delta^3}{a^4}\bigg)-d\bigg)\\
&\bigg(\beta\delta\gamma\bigg(\frac{1}{a}-\frac{\delta}{a^2}+\frac{\delta^2}{a^3}-\frac{\delta^3}{a^4}\bigg)-d\gamma\bigg)\frac{t^{2n+m}}{\Gamma(2n+m+1)},
\end{split}
\end{equation}
\begin{equation}\nonumber
\begin{split}
y_3(t)
=&\bigg(\beta\delta\gamma\bigg(\frac{2\delta}{a^3}-\frac{1}{a^2}-\frac{3\delta^2}{a^4}\bigg)+\beta\gamma\bigg(\frac{1}{a}-\frac{\delta}{a^2}+\frac{\delta^2}{a^3}-\frac{\delta^3}{a^4}\bigg)\bigg)\bigg(r-\frac{2r\delta}{K}-\alpha\delta\gamma\bigg(\frac{2\delta}{a^3}-\frac{1}{a^2}-\frac{3\delta^2}{a^4}\bigg)\\
&-\alpha\gamma\bigg(\frac{1}{a}-\frac{\delta}{a^2}+\frac{\delta^2}{a^3}-\frac{\delta^3}{a^4}\bigg)\bigg)\bigg(r\delta-\frac{r\delta^2}{K}-\alpha\delta\gamma\bigg(\frac{1}{a}-\frac{\delta}{a^2}+\frac{\delta^2}{a^3}-\frac{\delta^3}{a^4}\bigg)\bigg)\frac{t^{2m+n}}{\Gamma(2m+n+1)}\\
&-\alpha\delta\bigg(\frac{1}{a}-\frac{\delta}{a^2}+\frac{\delta^2}{a^3}-\frac{\delta^3}{a^4}\bigg)\bigg(\beta\delta\gamma\bigg(\frac{2\delta}{a^3}-\frac{1}{a^2}-\frac{3\delta^2}{a^4}\bigg)\\
&+\beta\gamma\bigg(\frac{1}{a}-\frac{\delta}{a^2}+\frac{\delta^2}{a^3}-\frac{\delta^3}{a^4}\bigg)\bigg)
\bigg(\beta\delta\gamma\bigg(\frac{1}{a}-\frac{\delta}{a^2}+\frac{\delta^2}{a^3}-\frac{\delta^3}{a^4}\bigg)-d\gamma\bigg)\frac{t^{m+2n}}{\Gamma(m+2n+1)}\\
&+\bigg(\frac{\beta\delta\gamma}{a^3}-\frac{3\beta\delta^2\gamma}{a^4}\bigg)\bigg(r\delta-\frac{r\delta^2}{K}-\alpha\delta\gamma\bigg(\frac{1}{a}-\frac{\delta}{a^2}+\frac{\delta^2}{a^3}-\frac{\delta^3}{a^4}\bigg)\bigg)^2\frac{\Gamma(2m+1)t^{2m+n}}{\Gamma(m+1)^2\Gamma(2m+n+1)}\\
&+\beta\delta\bigg(\frac{1}{a}-\frac{\delta}{a^2}+\frac{\delta^2}{a^3}-\frac{\delta^3}{a^4}\bigg)\bigg(\beta\delta\gamma\bigg(\frac{1}{a}-\frac{\delta}{a^2}+\frac{\delta^2}{a^3}-\frac{\delta^3}{a^4}\bigg)-d\gamma\bigg)\frac{t^{2n}}{\Gamma(2n+1)}\\
&+\beta\gamma\bigg(\frac{1}{a}-\frac{\delta}{a^2}+\frac{\delta^2}{a^3}-\frac{\delta^3}{a^4}\bigg)\bigg(r\delta-\frac{r\delta^2}{K}-\alpha\delta\gamma\bigg(\frac{1}{a}-\frac{\delta}{a^2}+\frac{\delta^2}{a^3}-\frac{\delta^3}{a^4}\bigg)\bigg)\frac{t^{m+n}}{\Gamma(m+n+1)}\\
&+\beta\bigg(\frac{1}{a}-\frac{\delta}{a^2}+\frac{\delta^2}{a^3}-\frac{\delta^3}{a^4}\bigg)\bigg(r\delta-\frac{r\delta^2}{K}-\alpha\delta\gamma\bigg(\frac{1}{a}-\frac{\delta}{a^2}+\frac{\delta^2}{a^3}-\frac{\delta^3}{a^4}\bigg)\bigg)\\
&\bigg(\beta\delta\gamma\bigg(\frac{1}{a}-\frac{\delta}{a^2}+\frac{\delta^2}{a^3}-\frac{\delta^3}{a^4}\bigg)-d\gamma\bigg)\frac{\Gamma(m+n+1)t^{m+2n}}{\Gamma(m+1)\Gamma(n+1)\Gamma(m+2n+1)}\\
&+\bigg(\beta\delta\bigg(\frac{1}{a}-\frac{\delta}{a^2}+\frac{\delta^2}{a^3}-\frac{\delta^3}{a^4}\bigg)-d\bigg)\bigg(\beta\gamma\bigg(\frac{1}{a}-\frac{\delta}{a^2}+\frac{\delta^2}{a^3}-\frac{\delta^3}{a^4}\bigg)\\
&+\beta\delta\gamma\bigg(\frac{2\delta}{a^3}-\frac{1}{a^2}-\frac{3\delta^2}{a^4}\bigg)\bigg)
\bigg(r\delta-\frac{r\delta^2}{K}-\alpha\delta\gamma\bigg(\frac{1}{a}-\frac{\delta}{a^2}+\frac{\delta^2}{a^3}-\frac{\delta^3}{a^4}\bigg)\bigg)\frac{t^{m+2n}}{\Gamma(m+2n+1)}\\
&+\bigg(\beta\delta\bigg(\frac{1}{a}-\frac{\delta}{a^2}+\frac{\delta^2}{a^3}-\frac{\delta^3}{a^4}\bigg)-d\bigg)^2\bigg(\beta\delta\gamma\bigg(\frac{1}{a}-\frac{\delta}{a^2}+\frac{\delta^2}{a^3}-\frac{\delta^3}{a^4}\bigg)-d\gamma\bigg)\frac{t^{3n}}{\Gamma(3n+1)},
\end{split}
\end{equation}
and so on.\\
Thus, the 3rd order approximate solution is obtained as
\begin{equation}\label{Rosenzweig-MacArthur predator-prey model_9}
x(t) = \sum_{n=0}^3 x_n(t) = x_0(t) + x_1(t) + x_2(t) + x_3(t), \\
\end{equation}
\begin{equation}\label{Rosenzweig-MacArthur predator-prey model_10}
y(t) = \sum_{n=0}^3 y_n(t) = y_0(t) + y_1(t) + y_2(t) + y_3(t), \\
\end{equation}
One can calculate more terms in a similar way to obtain better approximation of the solution.
\newpage
\section{Numerical computations}
\begin{figure}[H]
	%\centering
	\includegraphics[width=3.in, height=2.5in]{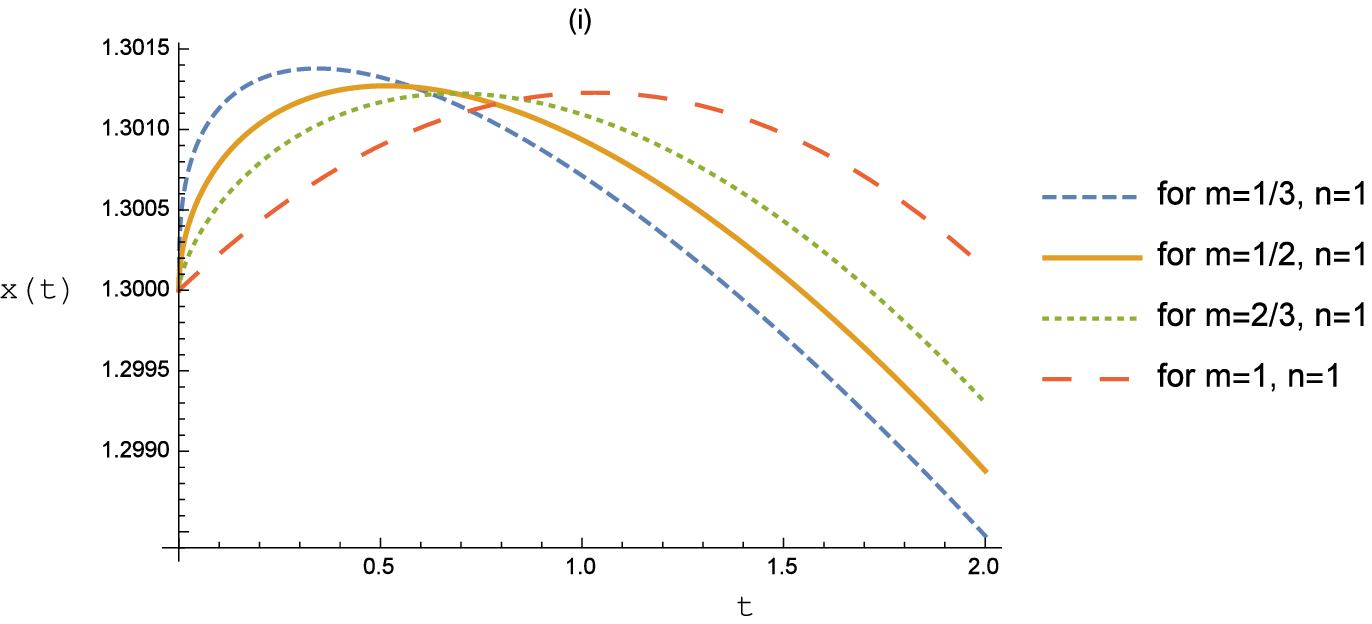}
	\includegraphics[width=3.in, height=2.5in]{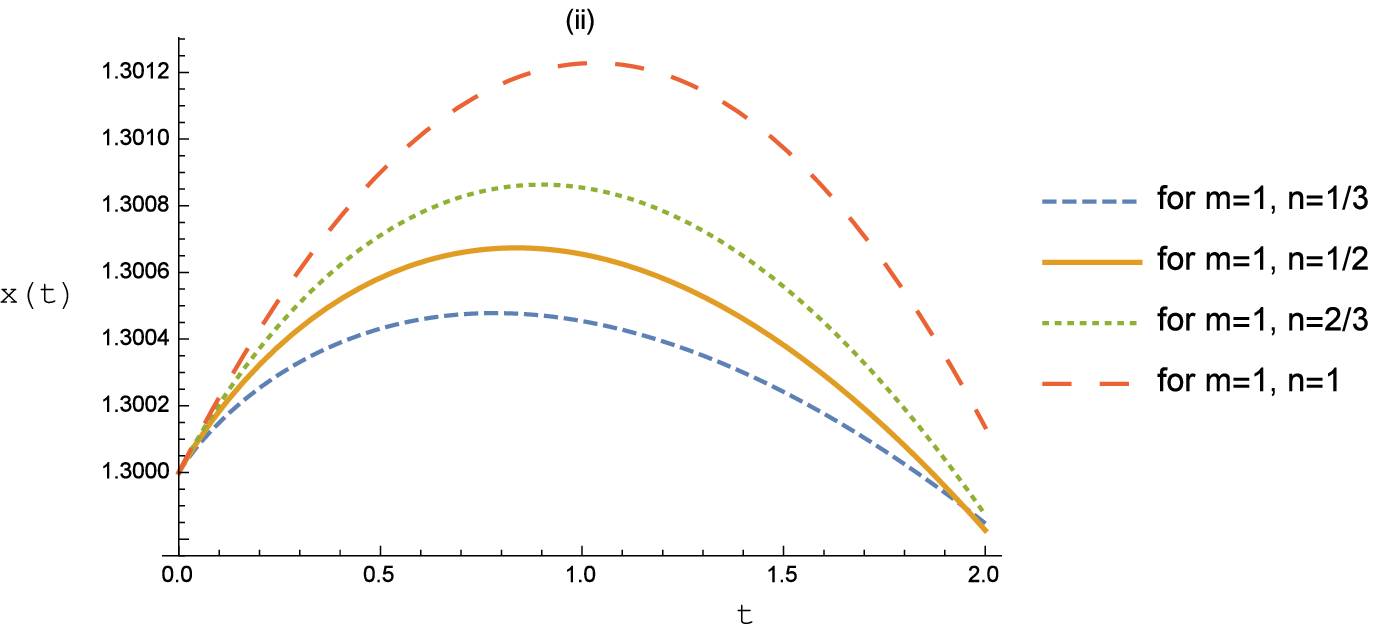}
	\caption{Approximate solutions of prey population, $x(t)$, of the system $(1.5)$ for some fixed fractional orders: (i) $m = 1/3, 1/2, 2/3, 1$ and $n = 1$ (ii) $m = 1$ and $n = 1/3, 1/2, 2/3, 1$.}
	\label{Prey.eps}
\end{figure}

\begin{figure}[H]
	%\centering
	\includegraphics[width=3in, height=2.5in]{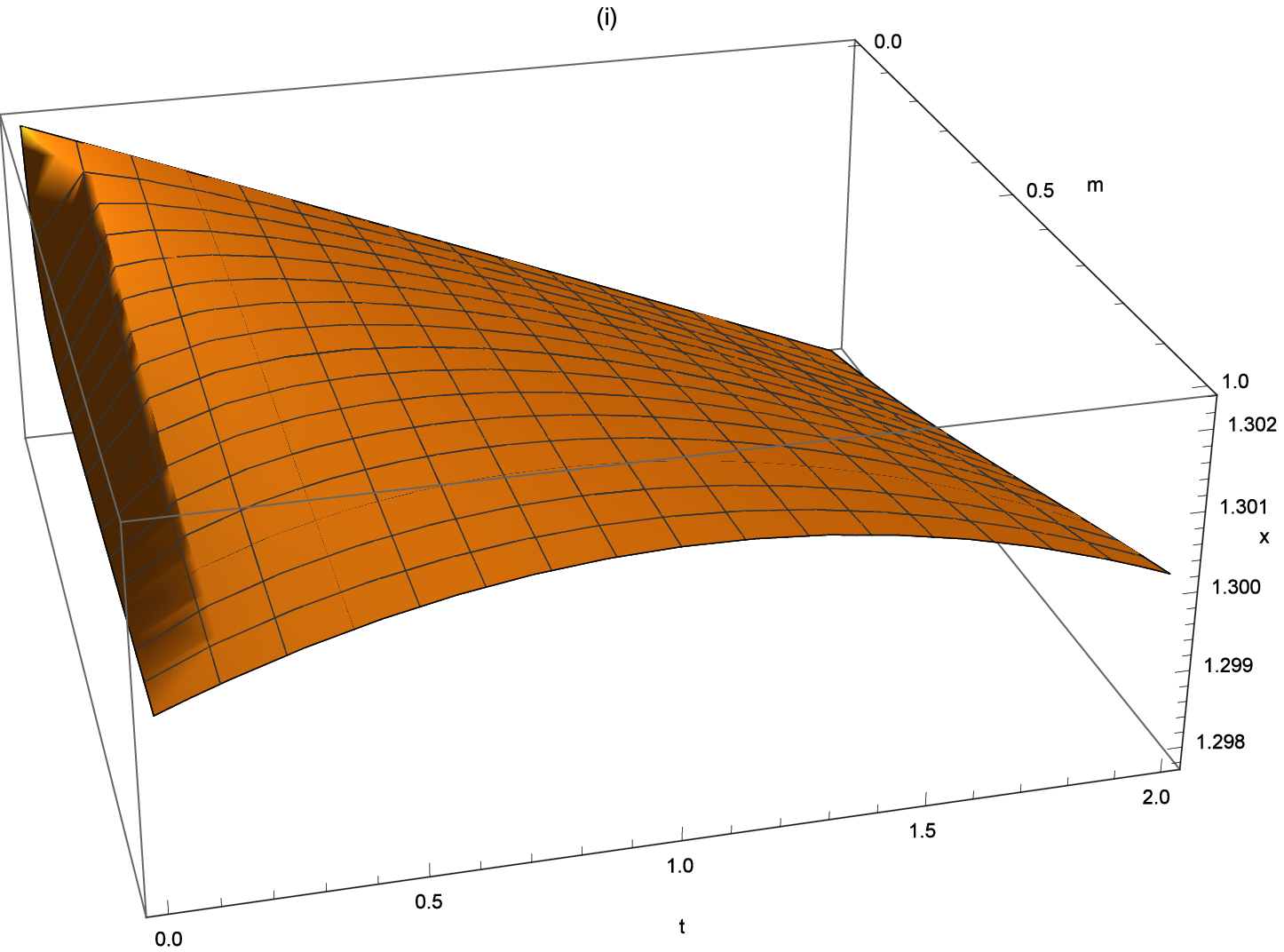}
	\includegraphics[width=3in, height=2.5in]{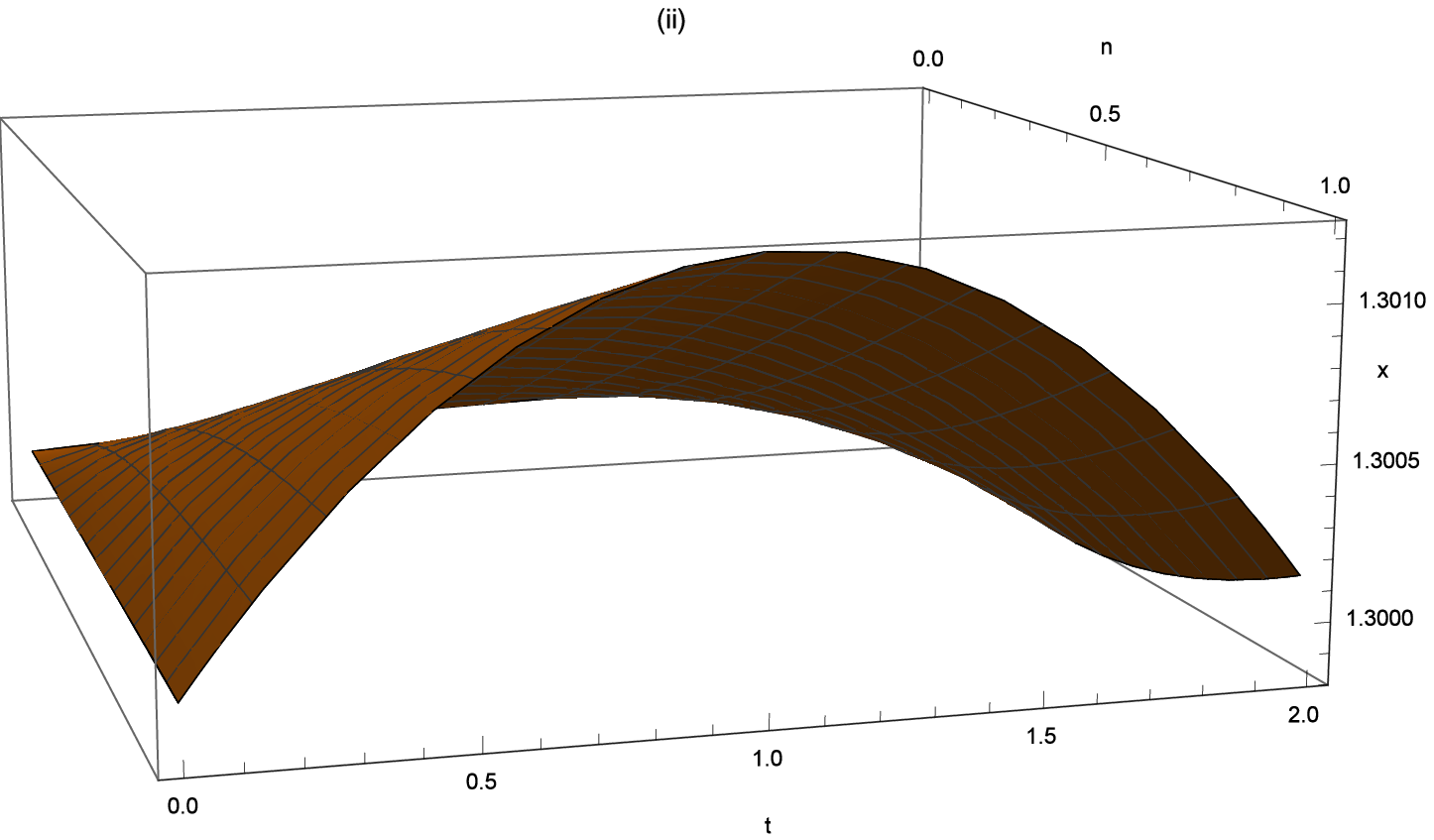}
	\caption{(i) Approximate solutions of prey population, $x(t)$, when both $m$ and $t$ vary with fixed value of $n = 1$. (ii) Approximate solutions of prey population, $x(t)$, when both $n$ and $t$ vary with fixed value of $m = 1$.}
	\label{Prey_3D.eps}
\end{figure}

In this section, we perform numerical computations of our model system (\ref{Rosenzweig-MacArthur predator-prey model}) for different fractional orders $(m,n = 1/3,1/2,2/3)$ and  as well as for the standard order $(m,n = 1)$.   We consider the parameter values as $r = 0.03, K = 10, a = 16, \alpha = 0.7, \beta = 0.6, d = 0.01$ and initial values as $\delta = 1.3$ and $\gamma = 0.6$. With this parameter set, we plot approximate solutions of prey population, $x(t)$, for different fractional orders $m$ when $n = 1$ (Fig. \ref{Prey.eps}(i)). This figure shows that prey population $x(t)$ reaches to its maximum more rapidly with decreasing m. Subsequent rapid decrement in population density is also observed in this case. Similar approximate solution of prey population $x(t)$  are plotted for different $n$ when $m = 1$ (Fig. \ref{Prey.eps}(ii)). This figure shows that maximum prey density increases with increasing $n$.
The corresponding surface plots are given in Fig. \ref{Prey_3D.eps}(i) when $t$ and $m$ varies with fixed value of $n$ $(n = 1)$ and in Fig. \ref{Prey_3D.eps}(ii) when $t$ and $n$ varies with fixed value of $m$ $(m = 1)$.

\begin{figure}[H]
	%\centering
	\includegraphics[width=3in, height=2.5in]{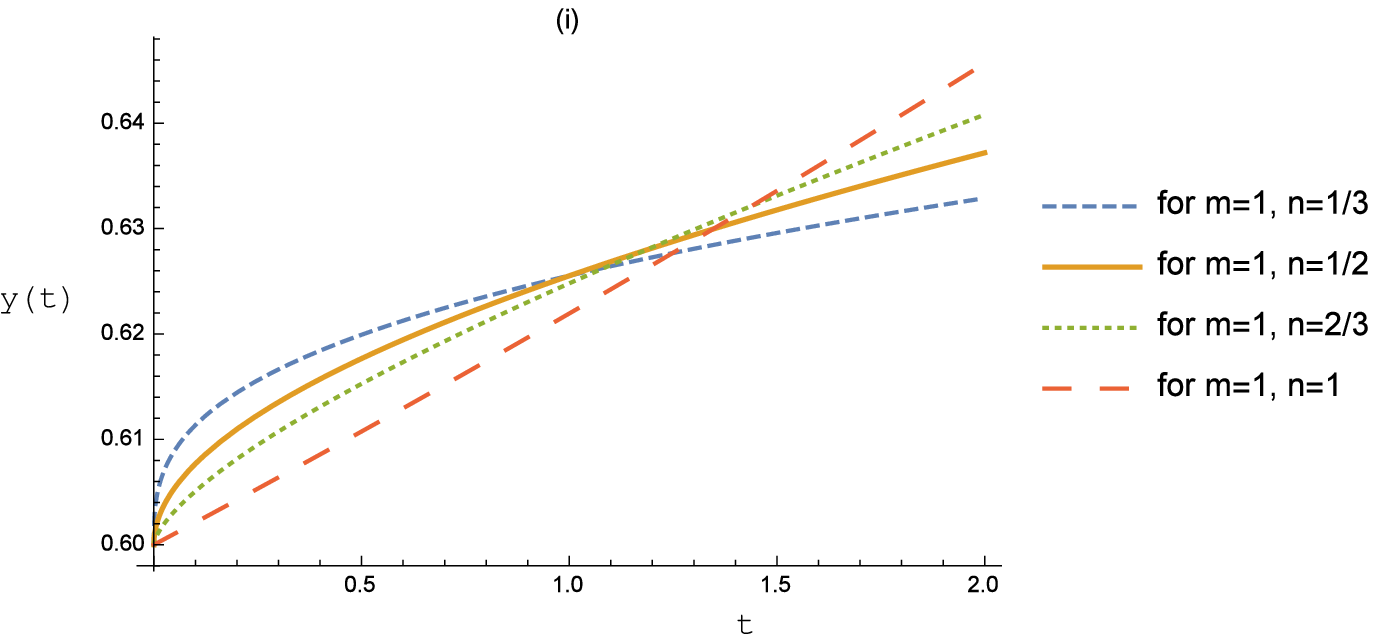}
	\includegraphics[width=3in, height=2.5in]{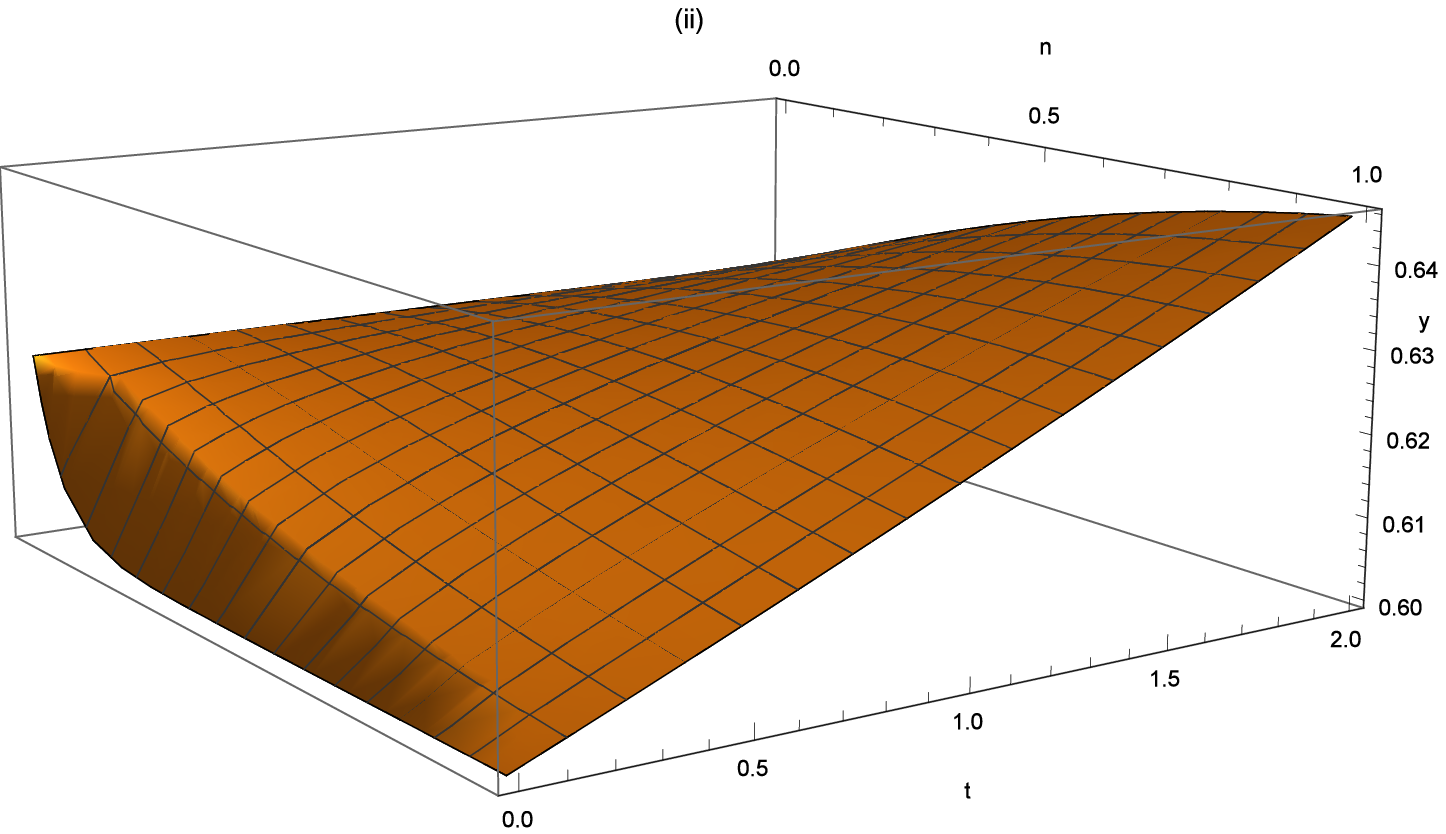}
	\caption{(i) Approximate solutions of predator population, $y(t)$, of the system $(1.5)$ for some fixed values of $n = 1/3, 1/2, 2/3, 1$ with $m = 1$. (ii) Approximate solutions of predator population, $y(t)$, when both $n$ and $t$ vary with fixed value of $m$ $(m = 1)$.}
	\label{Predator.eps}
\end{figure}
Approximate solutions of predator population, $y(t)$, for  different values of $n$ with $m = 1$ are presented in Fig. \ref{Predator.eps}(i). This figure shows that predator population grows very fast initially but satiates as time increases with decreasing values of $n$. The corresponding surface plot is presented in Fig. \ref{Predator.eps}(ii) when both $n$ and $t$ varies with fixed value of $m$ $(m = 1)$. One can draw similar graphs for predator population $y(t)$ when $m$ takes different values between $0$ and $1$ but $n$ is fixed at $1$ (figures not shown).\\
\begin{figure}[H]
	%\centering
	\includegraphics[width=3in, height=2in]{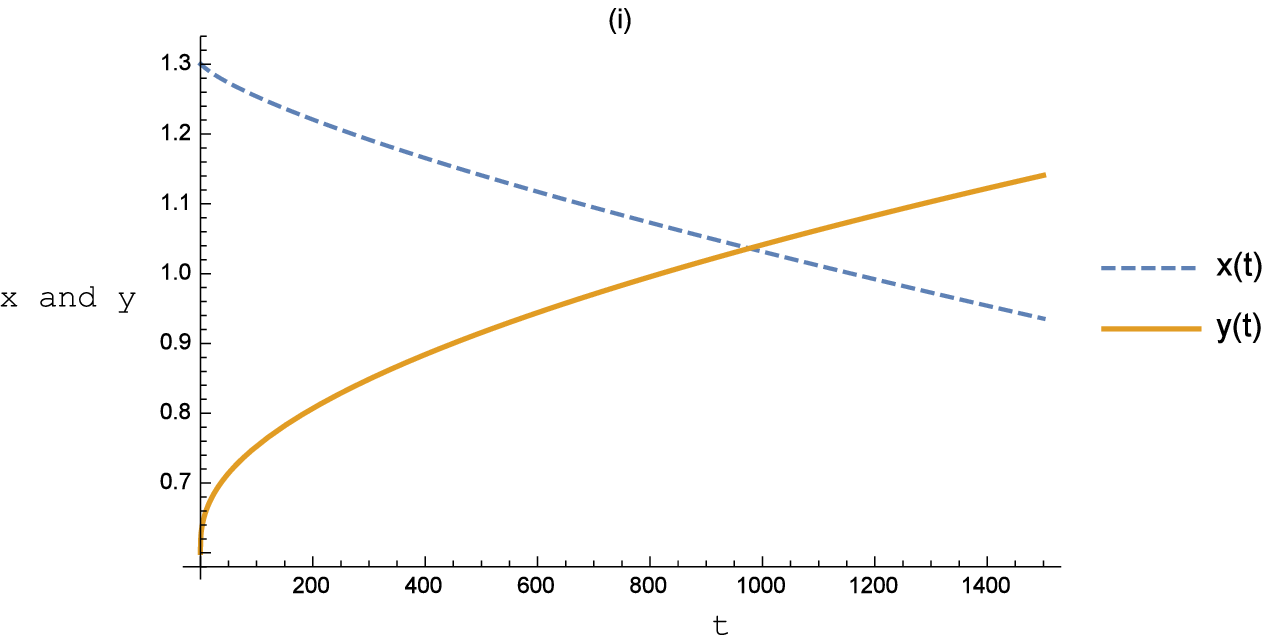}
	\includegraphics[width=3in, height=2in]{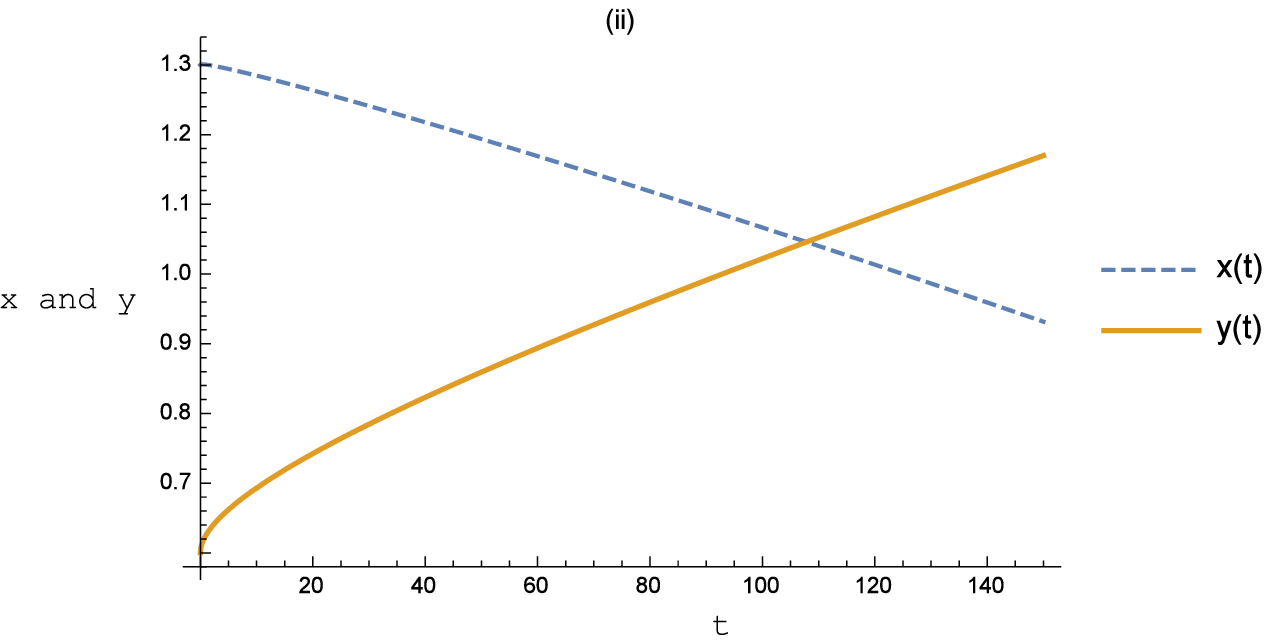}\\
	\includegraphics[width=3in, height=2in]{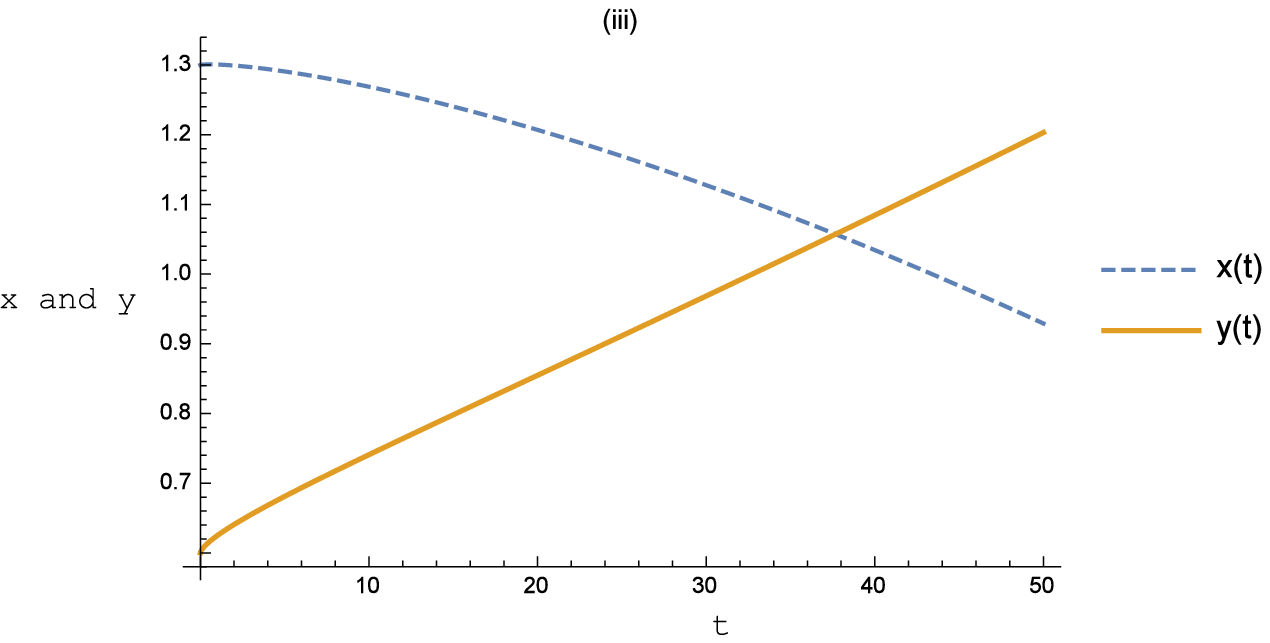}
	\includegraphics[width=3in, height=2in]{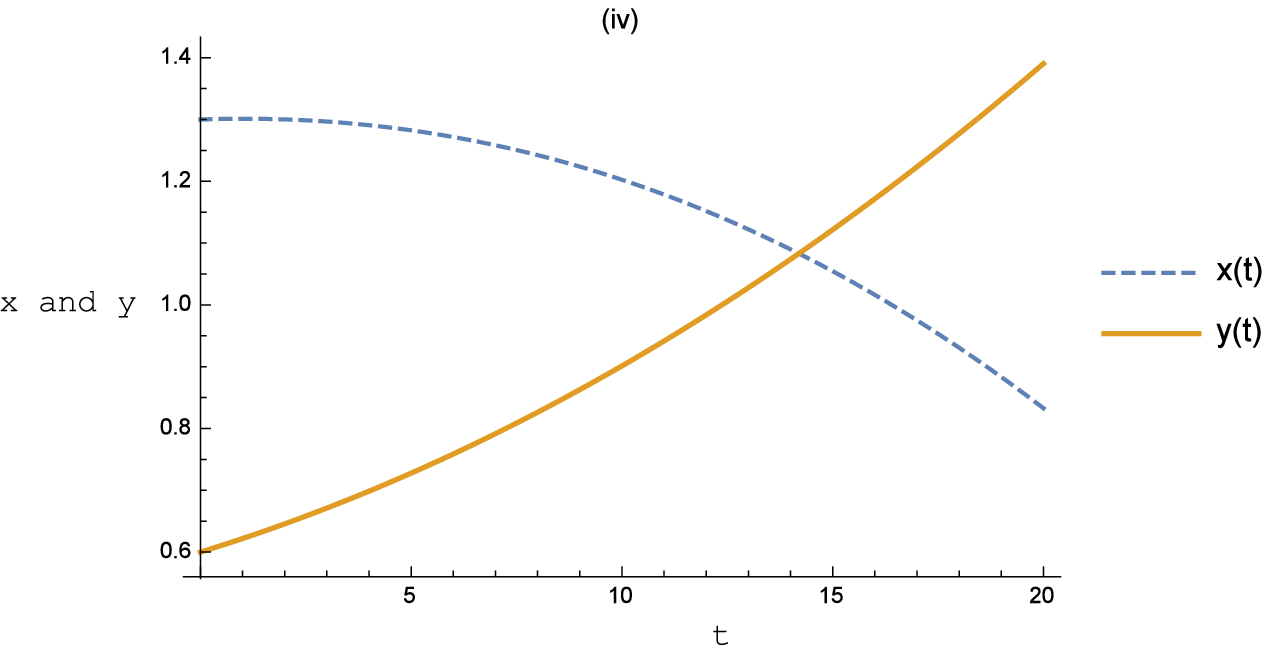}
	\caption{Approximate solutions of $x(t)$ and $y(t)$ for some fixed values of $m$ and $n$: (i)~$m = n = 1/3$, (ii)~$m = n = 1/2$, (iii)~$m = n = 2/3$, (iv)~$m = n = 1$.}
	\label{Prey_predator_dual.eps}
\end{figure}
In Fig. \ref{Prey_predator_dual.eps}, we simultaneously plot $x(t)$ and $y(t)$ for different orders. These figures show that prey population always decreases and predator population always increases whatever be the orders of $m$ and $n$. However, predator population grows more slowly and prey population declines more slowly as the orders of $m$ and $n$ increase.
\section{Summary}
Integer order differential equations are extensively used to study predator-prey model \cite{Freedman80}. However, application of fractional order differential equations is scarce. In recent time, researchers are showing interest to apply fractional time derivative in different field of natural sciences, including biology \cite{Das11}. In this article, we have used a system of fractional order differential equations to represent a predator-prey interaction. It is assumed that prey population grows logistically in absence of predator and predator consumes prey following satiated type II response function. Similar fractional order predator-prey models were studied by other researchers where predation term has been considered as unsaturated type \cite{CuiYang14,CuiYang12,DasGupta09,DasGupta11}, which is biologically unrealistic. To the best of our knowledge, no body has tried to find analytical solution of a predator-prey model where predator's functional response is saturated type II. We have used homotopy perturbation technique, which is supposed to be more efficient method compare to other methods like Adomian decomposition method (ADM) \cite{CuiYang14}, to explicitly find the analytical solution of the nonlinear predator-prey fractional order system in series form. It is our believe that the idea of solving a system of fractional order population model with type II response function will motivate other researchers to find analytical solutions of more complex and realistic biological systems.\\

%\noindent\textbf{References}\\


\begin{thebibliography}{0}
	
	\bibitem{CuiYang14} Z. Cui and Z. Yang, Homotopy perturbation method applied to the solution of fractional lotka-volterra equations with variable coefficients, Journal of Modern Methods in Numerical Mathematics 5 (2014), no. 1, 1-9.
	
	
	\bibitem{CuiYang12} Z. Yang Z. Cui and Z. Rui, Application of homotopy perturbation method to nonlinear fractional population dynamics models, Int. J. Appl. Math. Comput. 4 (2012), 403-412.
	
	
	
	\bibitem{May72} R. M. May, Limit cycles in prey predator communities, Science 177 (1972), 900-902.
	
	
	\bibitem{MomaniOdibat07} S. Momani and Z. Odibat, Homotopy perturbation method for nonlinear partial differential equations of fractional order, Physics Letters A 365 (2007), 345-350.
	
	
	\bibitem{Das11} S. Das, Introduction to fractional calculus for scientists and engineers, Springer, 2011.
	
	\bibitem{Podlubny99} I. Podlubny, Fractional differential equations, vol. 198, Academic Press, San Diego, Calif, USA, 1999
	
	
	\bibitem{DasGupta09} P.K. Gupta S. Das and Rajeev, A fractional predator prey model and its solution, Int. J. Nonlin. Sci. Numer. Simul. 10 (2009), 873-876.
	
	
	\bibitem{DasGupta11} S. Das, P.K. Gupta and Rajeev, A mathematical model on fractional lotka-volterra equations, J. Theoret. Bio. 277 (2011), 1-6.
	
	
	\bibitem{Freedman80} H.I. Freedman, Deterministic mathematical models in population ecology, Marcel Dekker, New York, 1980.
	
	\bibitem{HE99} J.H. He, Homotopy perturbation technique, Comput. Methods Appl. Mech. Engrg 178 (1999), 257.
	
	\bibitem{HE00} J.H. He, A coupling method of a homotopy technique and a perturbation technique for non linear problems, Int. J. Non-Linear Mech. 35 (2000), no. 1, 37.
	
	\bibitem{MillerRoss93} B. Ross K.S. Miller, An introduction to the fractional calculus and fractional differential equations, John Wiley and Sons, New York, 1993.
	
	\bibitem{Taghipour11} R. Taghipour, Application of homotopy perturbation method on some linear and nonlinear parabolic equations, IJRRAS 6 (2011), no. 1, 55-59.
	
	\bibitem{HE04} J. H. He, The homotopy perturbation method for nonlinear oscillators with discontinuities, Applied Mathematics and Computation 151 (2004), no. 1, 287-292. 
	
	\bibitem{HE05} J. H. He, Application of homotopy perturbation method to nonlinear wave equations, Chaos, Solitons and Fractals 26 (2005), no. 3, 695-700.
	
	
	\bibitem{HE98} J. H. He, An approximate solution technique depending on an artificial parameter: A special example, Communications in Nonlinear Science and Numerical Simulation 3 (1998), no. 2, 92-97.
	
	
	\bibitem{Ganji06} D. D. Ganji, The application of he’s homotopy perturbation method to nonlinear equations arising in heat transfer, Physics Letters A 353 (2006), no. 4-5, 337-341.
	
	\bibitem{BET07} J. Chattopadhyay N. Bairagi, P. K. Roy, Role of infection on the stability of a predator-prey system with several response functions - a comparative study, J. Theo. Biology 248 (2007), 10-25.
	
	\bibitem{SS12} Y. Satto J. Sugie, Uniqueness of limit cycles in a rosenzweig-mcarthur model with prey immigration, J. Theo. Biology 72 (2012), 299-316.
	\end{thebibliography}
\end{document}